\documentstyle[12pt]{article}
\newlength{\abstractwidth}
\renewcommand{\thefootnote}{\fnsymbol{footnote}}
\renewcommand{\thanks}[1]{\footnote{#1}} 
\newcommand{\starttext}{ \setcounter{footnote}{0}
\renewcommand{\thefootnote}{\arabic{footnote}}}

\newcommand{\be}{\begin{equation}}
\newcommand{\bea}{\begin{eqnarray}}
\newcommand{\eea}{\end{eqnarray}} \newcommand{\ee}{\end{equation}}
 
 \def\ba{\begin{eqnarray}}
\def\ea{\end{eqnarray}}

\def\o{\omega}

\def\tr{{\rm tr}}
\def\det{{\rm det}}

\def\log{\,{\rm log}\,}
\def\exp{\,{\rm exp}\,}

\def\o{\omega}

\def\o{\omega}

\def\ge{\geq}
\def\le{\leq}

\def\p{\partial}

\def\[{{\bf [}}
\def\]{{\bf ]}}

\def\ddbar{i\p\bar\p}

\def\ric{{\rm Ric}}
\def\mathbb{\bf}

\def\eqref{\ref}

\newcommand{\neweqref}[1]{(\ref{#1})}

\begin{document}
\starttext \baselineskip=18pt \setcounter{footnote}{0}
\newtheorem{theorem}{Theorem}
\newtheorem{lemma}{Lemma}
\newtheorem{corollary}{Corollary}
\newtheorem{definition}{Definition}
\newtheorem{conjecture}{Conjecture}
\newtheorem{proposition}{Proposition}


\begin{center}
{\Large \bf UNIFORM ENTROPY AND ENERGY BOUNDS FOR FULLY NON-LINEAR EQUATIONS
\footnote{Work supported in part by the National Science Foundation under grant DMS-22-03273.}}

\medskip
\centerline{Bin Guo and  Duong H. Phong}

\medskip

\begin{abstract}

{\footnotesize Energy bounds which are uniform in the background metric are obtained from upper bounds for entropy-like quantities. The argument is based on auxiliary Monge-Amp\`ere equations involving sublevel sets, and bypasses the Alexandrov-Bakelman-Pucci maximum principle. In particular, it implies uniform $L^\infty$ bounds for systems coupling a fully non-linear equation to its linearization, generalizing the cscK equation.}

\end{abstract}

\end{center}

\baselineskip=15pt
\setcounter{equation}{0}
\setcounter{footnote}{0}

\section{Introduction}
\setcounter{equation}{0}

Let $(X,\omega_X)$ be a compact $n$-dimensional K\"ahler manifold. 
If $\varphi$ is any smooth $\omega_X$-plurisubharmonic function, its entropy ${\rm Ent}(\varphi)$ is defined as the entropy of the measure $\o_\varphi^n=(\o_X+i\p\bar\p\varphi)^n$ with respect to the measure $\o_X^n$, 
\bea
{\rm Ent}(\varphi)=\int_X \log({\o_\varphi^n\over\o_X^n}) \o_\varphi^n
\eea
and, if $\varphi$ is normalized so that ${\rm sup}_X\varphi=0$, its energy $E(\varphi)$ is defined as
\bea
E(\varphi)=\int_X(-\varphi)\o_\varphi^n=\|\varphi\|_{L^1(\o_\varphi^n)}.
\eea
Both notions are essential for the study of the Monge-Amp\`ere equation and the problem of constant scalar curvature K\"ahler metrics. For example, the entropy is the leading term in the Mabuchi functional, while the energy is closely related to the well-known Aubin-Yau functionals $I(\varphi)$ and $J(\varphi)$ of K\"ahler geometry. In a recent major breakthrough, it had actually been shown by X.X. Chen and J.R. Cheng \cite{CC} that, in the constant scalar curvature equation, an upper bound for the entropy is equivalent to bounds for $\varphi$ to all orders. In another direction, it has been known for some time that bounds for the entropy imply bounds for the energy \cite{BBEGZ19}, and more precise embeddings of spaces of potentials with finite entropy into $L^p(\o_\varphi^n)$ have now been established in \cite{DGL}.

\smallskip
Our interest in entropy-like quantities comes from a related but different source: entropy-like quantities can also be defined for general fully non-linear equations on K\"ahler manifolds, and they turn out to be central to the existence of a priori $L^\infty$ estimates \cite{GPT}. More precisely, let $\o$ be a K\"ahler form on $X$, and consider an equation of the form
\bea
\label{eqn:main}
f(\lambda[h_\varphi])=c_\o e^{F_\o},
\quad {\rm sup}_X\varphi=0,
\quad
\lambda[h_\varphi]\in\Gamma,
\eea
where $(h_\varphi)^j{}_k=\o_X^{j\bar m}(\o_\varphi)_{\bar m k}$ is the relative endomorphism between $\o_X$ and $\o_\varphi=\o+i\p\bar\p\varphi$, $\lambda[h_\varphi]$ the un-ordered vector of eigenvalues of $h_\varphi$, and $f(\lambda)$ a given function defined on a cone $\Gamma\subset {\bf R}^n$ satisfying the conditions (1-4) spelled out in \S 2 below. The function $F_\o$ is normalized to satisfy $\int_Xe^{nF_\o}\o_X^n=\int_X\o_X^n$, and we set $V_\o=\int_X\o^n$ to be the volume of $\omega$.  
The case considered in \cite{GPT} is when $\o=\chi+t\o_X$, $t\in (0,1]$, where $\chi$ is a given non-negative closed $(1,1)$-form, and the estimates were required to be uniform in $t$. This case includes the one corresponding to a fixed background metric, which can be obtained by setting $\chi=0$ and $t=1$.
It is then proved in Theorem 2, \cite{GPT} that for such $\o$ and any $p>n$, we must have
\bea
{\rm sup}_X|\varphi|\leq C
\eea
where $C$ depends only on $\o_X,\chi, n,p,\gamma$, and upper bounds for the following three quantities
\bea
{c_\o^n\over V_\o}, \quad E(\o)={c_\o^n\over V_\o}\int_X (-\varphi)f(\lambda[h_\o])^n\o_X^n,\quad {\rm Ent}_p(e^{nF_\o})=\int_X e^{nF_\o}|F_\o|^p\o_X^n.
\eea
In the particular case of the Monge-Amp\`ere equation $f(\lambda)=(\prod_{j=1}^n\lambda_j)^{1\over n}$, $\Gamma=\{\lambda_j>0, 1\leq j\leq n\}$, the first two quantities in the above list can be bounded by elementary arguments, so we obtain uniform bounds for the complex Monge-Amp\`ere equation depending only on an upper bound for ${\rm Ent}_p(\o)$ for any fixed $p>n$, thus recovering the classic estimate of Kolodziej \cite{K}, as well as the uniform version established by Demailly and Pali \cite{DP} and Eyssidieux, Guedj, and Zeriahi \cite{EGZ}. 

\smallskip
The proof of the $L^\infty$ bound for fully non-linear PDE's which we just described was based on a comparison with an auxiliary Monge-Amp\`ere equation involving integrals $A_s$ on the sublevel sets $\Omega_s=\{z\in X;\varphi<-s\}$ of the function $\varphi$. This method turns out to be particularly effective: it has been extended by various authors to stability estimates \cite{GPT2}, nef classes \cite{GPTW1}, moduli of continuity \cite{GPTW2, GS}, lower bounds for the Green's function \cite{GPS, GPSS}, as well as Hermitian manifolds \cite{GP} and parabolic equations \cite{CC1}.

\smallskip
However, a natural question concerning the above general $L^\infty$ estimate was whether {\it uniform} bounds for the energy $E(\o)$ can be obtained from bounds for the entropy-like quantity ${\mathrm Ent}_p(\o)$.  For fixed background metric $\o$, this had been done in Theorem 3, \cite{GPT}, using an argument inspired by Chen-Cheng \cite{CC} which relied on the Alexandrov-Bakelman-Pucci (ABP) maximum principle. The ABP maximum principle is a powerful method pioneered by Blocki (\cite{B}; see more applications in e.g. \cite{P,Sz}), but its dependence on the background metric can be delicate, so uniform energy bounds are still lacking.

\smallskip

The main goal of the present paper is to supply such uniform bounds for the energy $E(\o)$. It turns out that that they can be obtained once again by a modification of the auxiliary Monge-Amp\`ere equation involving sublevel sets of $\varphi$ used in \cite{GPT}. This auxiliary Monge-Amp\`ere equation bypasses the ABP maximum principle, and does yield uniform estimates. 
As an indirect consequence, it can be used to simplify parts of the arguments in \cite{CC}, and generalize the $C^0$ estimates there to uniform estimates as well. We provide a description of the precise results in the next section.

\section{Statement of the main results}
\setcounter{equation}{0}

We begin by stating the precise conditions on the nonlinear operator $f(\lambda)$. As in \cite{GPT}, we require that $f:\Gamma\to {\bf R}_+$ satisfies

 (1) $\Gamma\subset {\mathbb R}^n$ is a symmetric cone with 
\begin{equation}\label{eqn:cone}
\Gamma_n\subset \Gamma \subset \Gamma_1;
\end{equation}
Here $\Gamma_k$ is the cone of vectors $\lambda$ with $\sigma_j(\lambda)>0$ for $1\leq j\leq k$, where $\sigma_j(\lambda)$ is the $j$-th symmetric polynomial in $\lambda$. In particular, $\Gamma_1$ is the half-space defined by $\lambda_1+\cdots+\lambda_n>0$, and $\Gamma_n$ is the first octant, defined by $\lambda_j>0$ for $1\leq j\leq n$.

(2) $f(\lambda)$ is symmetric in $\lambda = (\lambda_1,\ldots, \lambda_n)\in \Gamma$ and it is homogeneous of degree one;

 (3) $\frac{\partial f}{\partial \lambda_j}>0$ for each $j=1,\ldots, n$ and $\lambda\in \Gamma$;

 (4) There is a $\gamma>0$ such that 
\begin{equation}\label{eqn:structure}
\prod_{j=1}^n \frac{\partial f(\lambda)}{\partial \lambda_j}\ge \gamma,\quad \forall \lambda\in \Gamma.
\end{equation}

It is well-known that equations such as the Monge-Amp\`ere equation, with $f(\lambda)=(\prod_{j=1}^n\lambda_j)^{1\over n}$, or the Hessian equation with $f(\lambda)=\sigma_k(\lambda)^{1\over k}$,
or the $p$-Monge-Amp\`ere equation of Harvey and Lawson \cite{HL1,HL2} with
$$f(\lambda) = \Big( \prod_{I} \lambda_I\Big)^{\frac{n!}{(n-p)!p!}}$$
 where $I$ runs over all distinct multi-indices $1\le {i_1}<\cdots < {i_p}\le n$, $\lambda_I = \lambda_{i_1} + \cdots + \lambda_{i_p}$, and $\Gamma$ is the cone defined by $\lambda_I>0$ for all  $p$-indices $I$, all satisfy the structural condition (4). In a remarkable recent development, Harvey and Lawson \cite{HL3} showed that the condition (4) actually holds for very large classes of non-linear operators, including all invariant Garding-Dirichlet operators. As noted in \cite{HL3}, the condition (4) also arose independently in \cite{AO} in the study of $W^{2,p}$ interior regularity.

\begin{theorem}\label{thm:main1}\label{thm:main}
Let $(X,\o_X)$ is a compact $n$-dimensional K\"ahler manifold. Let $\o$ be any K\"ahler form on $X$ with
\bea
\label{kappa}\label{eqn:bound}
\o\leq \kappa\,\o_X
\eea
for some constant $\kappa>0$. Consider the equation \neweqref{eqn:main} with the operator $f(\lambda)$ satisfying the conditions (1-4).
Then for any $p>0$, any $C^2$ solution $\varphi$ of \neweqref{eqn:main} satisfies the following 

{\rm (i)} Trudinger-like inequalities 
\begin{equation}\label{eqn:inequality}
\int_X e^{\alpha (-\varphi)^{q}} \omega_X^n\le C_T,
\end{equation} 

{\rm (ii)} and energy-like inequalities
\begin{equation}\label{eqn:inequality 2}
\int_X (-\varphi)^{pq} e^{nF_\omega} \omega_X^n\le C_e.
\end{equation}

Here the exponent $q$ is given by $q = \frac{n}{n-p}$ if $p<n$, and can be any strictly positive exponent if $p \geq n$. The constants $C_T$ and $C_e$ are computable constants depending only on $n, p, q, \omega_X, \kappa,\gamma$, and upper bounds for the following two quantities
\bea
{c_\o^n\over V_\omega},
\quad
{\mathrm{Ent}}_p(e^{nF_\o}) = \int_X e^{nF_\o} |F_\o|^p\omega_X^n,
\eea
and the $\alpha>0$ in \neweqref{eqn:inequality} is a constant that depends only on $n, p, \gamma, {c^n_\o\over V_\omega}$ and $\kappa$.

\end{theorem}

\medskip
We observe that, in the case of a fixed background K\"ahler metric $\o$, this theorem was proved as Theorem 3 in \cite{GPT}. The point of the new theorem is to have uniform estimates, even as the background metric $\o$ may degenerate to the boundary of the K\"ahler cone. For this same reason, we include the case $p>n$ in the statement. When $p>n$ and the background K\"ahler form $\o$ is fixed, it follows from Theorem 1, \cite{GPT}, that the solution $\varphi$ of the equation is actually bounded, and the above Trudinger-like and energy-like inequalities follow at once. But here again, the existing results do not give the inequalities uniform in $\o$ that we seek.

\smallskip
To obtain estimates which are uniform with respect to the background metric $\o$, we have to improve on the proof for fixed $\o$ in \cite{GPT}, which was modeled on the arguments of \cite{CC} for the constant scalar curvature equation, and made essential use of the ABP maximum principle. What appears needed for uniform estimates is rather arguments in the spirit of Theorem 2, \cite{GPT}, and indeed, it turns out that these arguments can be adapted to the case at hand.

\smallskip
Theorem \ref{thm:main1} readily combines with Theorem 2, \cite{GPT} to give the following improvement, which we state for easy reference in the future:

\begin{theorem}
\label{thm:main2}
Let $(X,\o_X)$ be a compact $n$-dimensional K\"ahler manifold. Let $\o$ be a K\"ahler form satisfying the condition (\ref{kappa}) for a fixed constant $\kappa>0$, and consider the equation 
 \neweqref{eqn:main} with the operator $f(\lambda)$ satisfying the conditions (1-4). Then for any $p>n$, a $C^2$ solution $\varphi$ of the equation \neweqref{eqn:main} must satisfy
 \bea
 {\rm sup}_X|\varphi|\leq C
 \eea
 where $C$ is a constant depending only on $\o_X,n,p,\gamma,\kappa$, and upper bounds for the following two quantities
 \bea
 {c_\o^n\over V_\o}, \quad {\rm Ent}_p(e^{nF_\o}).
 \eea
\end{theorem}

We observe that in \cite{GPT}, Theorem 2 was stated for background K\"ahler metrics of the form $\o=\chi+t\o_X$, $t\in (0,1]$. However, as noted in \cite{GPSS}, the proof applies uniformly for background K\"ahler forms $\o$ satisfying (\ref{kappa}), as long as we allow a dependence of all relevant constants on the bound $\kappa$.

\smallskip
We would also like to note that, for the specific case of the Monge-Amp\`ere equation on strongly pseudoconvex domains in ${\bf C}^n$, a proof of $L^\infty$ estimates using the Monge-Amp`ere energy and corresponding Sobolev inequalities has been given in \cite{WWZ, WWZ1}.

\medskip

As we have stressed above, the key to the proof of Theorem \ref{thm:main1} is an argument bypassing the use of the ABP maximum principle. As such, it can also apply and simplify several of those parts in the paper of Chen-Cheng \cite{CC} which relied on the ABP maximum principle. As an illustration, we state here a $C^0$ estimate for a coupled system, generalizing the coupled system corresponding to the constant scalar curvature equation, which is uniform with respect to the background K\"ahler form $\o$.

\medskip
Let $(X,\o_X)$ be a compact $n$-dimensional K\"ahler manifold as before, $\o$ be a K\"ahler form, and let $\theta$ be a given smooth $(1,1)$-form on $X$. We consider the coupled system
\bea
\label{eqn:3.1}
&&
f(\lambda[h_\varphi]) = c_\omega\, e^{F_\o}, \quad  \sup_X \varphi = 0\nonumber\\
&&
 \Box_{\omega_\varphi} F_\o = - c_\theta + G^{i\bar j} \theta_{i\bar j}, 
\eea
where $G^{i\bar j} = \frac{\partial}{\partial h_{i\bar j}} \log f(\lambda[h_\varphi])$ is the linearized operator of $\log f(\lambda[h])$ which is positive definite by condition (3) in the definition of $f(\lambda)$, and  $\Box_{\omega_\varphi} F_\o = G^{i\bar j} (F_\o)_{i\bar j}$. The function $F_\o$ is again normalized (for the purpose of determining $c_\o$) to satisfy $\int_X e^{nF_\o} \omega_X^n = \int_X \omega_X^n$, and $c_\theta$ is a constant determined by the equation. We have then

\begin{theorem}
\label{thm:main3}
Assume that the function $f(\lambda)$ satisfies the conditions (1-4) spelled out at the beginning of this section,  and $\o$ satisfies the condition $\o\leq \kappa\,\o_X$ for some constant $\kappa$. 
Fix a number $p\in (0,n]$. We assume that 
\begin{equation}\label{eqn:assumption}
{\mathrm{Ent}}_p(e^{nF_\o}) = \int_X  |F_\o|^p e^{nF_\o}\omega_X^n\le K_1,\quad {\rm and}\quad \theta \ge - K_2 \omega
\end{equation}
for some constants $K_1>0$ and $K_2>0$. Then
\bea\nonumber
{\rm sup}_X|\varphi|\leq C,\, \mbox{and } \sup\nolimits_X F_\o \le C,
\eea
for a constant $C$ depending only on $\o_X,p,n,\gamma, \kappa, c_\omega^n/V_\omega, c_\theta,K_1$ and $K_2$. If we assume further that
\bea\label{eqn:K3}
\theta\leq K_3\o
\eea
then the function $F_\o$ is bounded from below by another constant depending further on $K_3$.

\end{theorem}

We remark that when $f(\lambda[h_\varphi]) = (\frac{\omega_\varphi^n}{\omega_X^n})^{1/n}$ and $\theta = \ric(\omega_X)$ is the Ricci curvature of $\omega_X$, the coupled system \neweqref{eqn:3.1} is the constant scalar curvature K\"ahler equation (cscK) studied in \cite{CC}. In this case, the constants $c_\o = V_\omega^{1/n}$, and $c_\theta\in {\mathbb R}$ depend only on the cohomology classes $[\omega]$ and $c_1(X)$. The lower bound for $F_\o$ was established in \cite{He}. One of the main results of \cite{CC} is that, assuming an upper bound for the entropy ${\mathrm Ent}_{p=1}(e^{nF})$, one can obtain a priori estimates for $\varphi$ of all orders.
What our result shows is that the $C^0$ bounds for $\varphi$ and $F_\o$  for this particular coupled system still hold even when $\omega$ degenerates to the boundary of the K\"ahler cone. We remark that the condition $-K_2\omega\le \theta\le K_3\omega$ in \neweqref{eqn:assumption} and \neweqref{eqn:K3} is not very restrictive for a degenerating family. For example, it holds for $\omega = \chi + t \omega_X$ and $\theta = - \chi$ for some nonnegative $(1,1)$-form $\chi$.

\medskip

\section{Proof of Theorem \ref{thm:main}}\label{section 2}
\setcounter{equation}{0}

For notational simplicity, we will omit the subscript $\omega$ in $F_\omega$ and simply write $F$ in this section. 
Suppose $\varphi\in C^2(X)$ solves the equation \neweqref{eqn:main} with $\sup_X\varphi = 0$. Since $\lambda[h_\varphi] \in \Gamma\subset \Gamma_1$, we have 
$$\tr_{\omega_X} \omega + \Delta_{\omega_X}\varphi >0.$$
The assumption \neweqref{eqn:bound} implies that $\Delta_{\omega_X} \varphi\ge -\tr_{\omega_X} \omega \ge - n\kappa$. An application of the Green's formula shows the uniform $L^1(X,\omega_X^n)$ estimate of $\varphi$, i.e.
\begin{lemma}
There exists a constant $C_0 = C_0(n,\kappa, \omega_X)$ such that 
$$\int_X |\varphi| \omega^n_X \le C_0.$$
\end{lemma}
We will write $K>0$ to be an upper bound of the $p$-th entropy of $e^{nF}$ for $p\in (0,n]$. As in Theorem \ref{thm:main}, let $q = \frac{n}{n-p}$ if $p<n$ and $q>0$ be any positive number if $p = n$. Let $s>0$ be a positive number and $\Omega_s\subset X$ be the sub-level set
\begin{equation}\label{eqn:level set}
\Omega_s = \{ z\in X| ~- \varphi(z) - s >0\}.
\end{equation}
We also define a monotonically decreasing function $\phi(s)$ as in \cite{GPT}
\begin{equation}\label{eqn:phi}
\phi(s) = \int_{\Omega_s} e^{nF} \omega_X^n.
\end{equation} 
Given these definitions, we have the following lemma about the decay of $\phi(s)$. 
\begin{lemma}\label{lemma 2}
There exists a constant $C_1>0$ that depends on $n, p, K, C_0>0$ such that for any $s>1$
\begin{equation}\label{eqn:lemma 2}
\phi(s)\le \frac{C_1}{(\log\,s)^p}.
\end{equation}
\end{lemma}
\noindent{\em Proof of Lemma \ref{lemma 2}.} We observe that by the H\"older-Young's inequality (c.f. \cite{GPT}), there is a constant $C_p>0$ depending only $p$ such that
\bea\nonumber
\int_{\Omega_s} e^{nF} \omega_X^n &\le & 2^p \int_{\Omega_s} \Big(\frac{2^{-1} \log(-\varphi)}{ \log s} \Big)^p e^{nF} \omega_X^n\\
&\le & \frac{2^p}{(\log s)^p} \int_{\Omega_s} \Big( e^{nF} (1+ n^p|F|^p) + C_p e^{\log (-\varphi)}     \Big) \omega_X^n\nonumber\\
&\le \nonumber& \frac{C}{(\log s)^p} ( 1+ K + C_p C_0  )\le \frac{C_1}{ (\log s)^p}.
\eea 
From Lemma \ref{lemma 2}, we see that $\phi(s)$ can be arbitrarily small if $s>1$ is sufficiently large. We now prove Theorem \ref{thm:main}.

\medskip

\noindent{\em Proof of Theorem \ref{thm:main}.} We will modify the approach given in \cite{GPT}. We break the proof into four steps.

\medskip

\noindent{\bf Step 1.} Let $\tau_k(x): {\mathbb R}\to {\mathbb R}_{+}$ be a sequence of positive smooth functions that converges monotonically decreasingly to the function $x\cdot\chi_{{\mathbb R}_{+}}(x)$. Let $a = pq = \frac{np}{n-p}$ ($a$ is {\em any} positive number if $p=n$). We solve the following complex Monge-Amp\`ere equation
\begin{equation}\label{eqn:aux}
(\omega + \ddbar \psi_{s,k })^n = \frac{\tau_k(-\varphi - s) ^a}{A_{s,k}} c_\omega^n e^{nF} \omega_X^n,\quad \sup_X \psi_{s,k} = 0.
\end{equation}
Here the constant $A_{s,k}$ is defined by 
\begin{equation}\label{eqn:Ask}
A_{s,k} = \frac{c_\omega^n}{V_\omega} \int_X \tau_k(-\varphi - s) ^ a e^{nF} \omega_X^n
\end{equation}
to make the equation \neweqref{eqn:aux} compatible.  By assumption $[\omega]$ is a K\"ahler class, so by Yau's theorem \cite{Y}, equation \neweqref{eqn:aux} admits a unique smooth solution $\psi_{s,k}$. We observe that by dominated convergence theorem that as $k\to\infty$
\begin{equation}\label{eqn:As}
A_{s,k}\to A_s: = \frac{c_\omega^n}{V_\omega} \int_{\Omega_s} (-\varphi - s)^a e^{nF} \omega_X^n. 
\end{equation}

\smallskip

\noindent{\bf Step 2.} Define a smooth function
$$\Phi := -\varepsilon (-\psi_{s,k} + \Lambda)^b - \varphi - s, $$
where the constants are given by \begin{equation}\label{eqn:constant 1}
b = \frac{n}{n+a}\in(0,1), \quad \mbox{and }\varepsilon = \frac{1}{\gamma^{1/(n+a)} (n b)^{n/(n+a)}} A_{s,k}^{\frac{1}{n+a}},
\end{equation}
and $\Lambda$ is chosen so that $\varepsilon b \Lambda^{-(1-b)} = 1$, i.e.
\begin{equation}\label{eqn:constant 2}
\Lambda = \frac{b^{1/(1-b)}}{(\gamma^{1/(n+a)} (nb)^{n/(n+a)})^{1/(1-b)}} A_{s,k}^{\frac{1}{a}}.
\end{equation}
We {\bf claim} that $\Phi\le 0$ on $X$. To see this, we note that $\Phi<0$ on $X\backslash \Omega_s$ by definition. If $\max_X \Phi$ is achieved somewhere on $X\backslash \Omega_s$, we are done. So we assume $\max_X \Phi = \Phi(x_0)$ for some point $x_0\in\Omega_s$. Let 
\begin{equation}\label{eqn:G}G^{i\bar j} = \frac{\partial \log f (\lambda[h_\varphi])}{\partial h_{i\bar j}} = \frac{1}{f} \frac{\partial f (\lambda[h_\varphi])}{\partial h_{i\bar j}}\end{equation}
be the coefficients of the linearized operator of $f(\lambda[h_\varphi])$ at the point $x_0$. $G^{i\bar j}$ is positive definite by condition (3) of the function $f$. Hence we have $G^{i\bar j} \Phi_{i\bar j}\le 0$ at the point $x_0$. Choosing local holomorphic coordinates at $x_0$ such that at this point $(\omega_X)_{i\bar j}|_{x_0} = \delta_{ij}$
 and $(\omega_\varphi)|_{x_0}$ is diagonal with eigenvalues $\lambda_1,\ldots, \lambda_n$. Then we have
 \begin{equation}\label{eqn:2.9}
 \det G^{i\bar j} = \frac{1}{f^n} \prod_{j=1}^n \frac{\partial f}{\partial \lambda_j} \ge \frac{\gamma}{f^n},
 \end{equation}
 by the condition (4) in the definition of the nonlinear operator $f(\lambda)$.
 
   We calculate at $x_0$ as follows:
 \bea\noindent
 0 &\ge &\nonumber G^{i\bar j}\Phi_{i\bar j} \\
 & = \nonumber& \varepsilon b (-\psi_{s,k} + \Lambda)^{b- 1} G^{i\bar j} (\psi_{s,k})_{i\bar j} + \varepsilon b (1-b) (-\psi_{s,k}+\Lambda)^{b-2} G^{i\bar j} (\psi_{s,k})_i (\psi_{s,k})_{\bar j}\\
 && \nonumber - G^{i\bar j} (\omega_\varphi - \omega)_{i\bar j}\\
 &\ge &  \nonumber \varepsilon b (-\psi_{s,k} + \Lambda)^{b- 1} G^{i\bar j} (\omega_{\psi_{s,k}})_{i\bar j} -1 + (1-  \varepsilon b (-\psi_{s,k} + \Lambda)^{b- 1}    ) G^{i\bar j} \omega_{i\bar j}\\
 &\ge & n \varepsilon b (-\psi_{s,k} + \Lambda)^{b- 1}  (\det G^{i\bar j})^{1/n} (\det (\omega_{\psi_{s,k}})_{i\bar j})^{1/n } - 1 + (1- \varepsilon b \Lambda^{-(1-b)}) G^{i\bar j} \omega_{i\bar j}\nonumber\\
 &\ge \nonumber& n \varepsilon b(-\psi_{s,k} + \Lambda)^{b-1} \gamma^{\frac 1 n}\Big( \frac{(-\varphi - s)^a}{A_{s,k}}     \Big)^{\frac 1 n} - 1.
 \eea
Here we have applied the arithmetic-geometric inequality, \neweqref{eqn:2.9}, and the equation \neweqref{eqn:aux} of $\omega_{\psi_{s,k}}$. It follows easily from the above and the choice of constants in \neweqref{eqn:constant 1} that at $x_0$
$$(-\varphi - s) \le \frac{A_{s,k}^{1/a}}{(  \gamma^{1/n} n \varepsilon b  )^{n/a}} (-\psi_{s,k} + \Lambda)^{(1-b) \frac{n}{a}} = \varepsilon (-\psi_{s,k} + \Lambda)^{b}$$
that is, $\Phi(x_0)\le 0$.

\medskip

\noindent {\bf Step 3.} From the previous step, $\Phi\le 0$. Thus on $\Omega_s$ we have 
\begin{equation}\label{eqn:2.10}
\frac{(-\varphi - s)}{A_{s,k}^{1/(n+a)}} \le C ( -\psi_{s,k} + C A_{s,k}^{1/a} )^{\frac{n}{n+a}},
\end{equation}
for some uniform constant $C>0$ that depends on $n, p,\gamma$. Taking $((n+a)p/n)$-th power on both sides of \neweqref{eqn:2.10}, and multiplying $e^{nF}$, we obtain that on $\Omega_s$
\begin{equation}\label{eqn:2.11}
\frac{(-\varphi - s)^{\frac{p(n+a)}{n}}}{A_{s,k}^{p/n}} e^{nF} \le C (-\psi_{s,k} + A_{s,k}^{1/a})^p e^{nF}\le C_2[ (-\psi_{s,k})^p e^{nF} + A_{s,k}^{p/a} e^{nF}  ],
\end{equation}
for some constant $C_2>0$ depending only on $n, p,\gamma$. We note that by H\"older-Young's inequality, for any $\beta>0$ there is a constant $C_p>0$ depending only on $p$ such that  
\begin{equation}\label{eqn:HY}
(-\frac {\beta }{2}\psi_{s,k})^p e^{nF} \le e^{nF} (1+ |nF|^p) + C_p e^{-\beta \psi_{s,k}}.
\end{equation}
Since $\omega+ \ddbar \psi_{s,k}>0$ and by the assumption \neweqref{eqn:bound}, that is $\omega\le \kappa \omega_X$, we see that $\psi_{s,k}\in PSH(X,\kappa \omega_X)$. Hence there exists a $\beta = \beta(X,\kappa \omega_X)>0$ such that (\cite{H, Ti})
\begin{equation}\label{eqn:2.13}\int_X e^{-\beta \psi_{s,k}} \omega_X^n \le C_X,\end{equation}
for some uniform constant $C_X = C_X(\kappa \omega_X, n)$. We integrate both sides of \neweqref{eqn:2.11} against $\omega_X^n$ over $\Omega_s$ and apply \neweqref{eqn:HY} and \neweqref{eqn:2.13},
\begin{equation}\label{eqn:2.14}
\int_{\Omega_s} \frac{(-\varphi - s)^{\frac{p(n+a)}{n}}}{A_{s,k}^{p/n}} e^{nF} \omega_X^n \le C_3 + C_2 A_{s,k}^{p/a} \int_{\Omega_s} e^{nF} \omega_X^n = C_3 + C_2 A_{s,k}^{p/a} \phi(s).
\end{equation}
Here $C_3>0$ is a uniform constant depending on $n, p, \kappa, C_X$ and $K$, the upper bound of ${\mathrm{Ent}}_p(e^{nF})$. Letting $k\to \infty$ in \neweqref{eqn:2.14}, we obtain that
\begin{equation}\label{eqn:2.15}
\int_{\Omega_s} (-\varphi - s)^{\frac{p(n+a)}{n}} e^{nF} \omega_X^n \le C_3 A_s ^{p/n}+ C_2 A_{s}^{\frac{p}{a} + \frac{p}{n}} \phi(s),
\end{equation}
where $A_s$ is given by \neweqref{eqn:As}. On the other hand, by H\"older inequality we have
\bea\nonumber
A_s&  = & \frac{c_\o^n}{V_\omega} \int_{\Omega_s} (-\varphi - s)^a e^{nF} \omega_X^n
\\
&\le &\nonumber \frac{c_\omega^n}{V_\omega} \Big( \int_{\Omega_s} (-\varphi - s)^{\frac{p(n+a)}{n}} e^{nF} \omega_X^n \Big) ^{\frac{na}{p(n+a)}} \Big( \int_{\Omega_s} e^{nF} \omega_X^n\Big) ^{1- \frac{na}{(n+a) p}}\\
&\le \nonumber &\frac{c_\omega^n}{V_\omega} \Big(C_3 A_s ^{p/n}+ C_2 A_{s}^{\frac{p}{a} + \frac{p}{n}} \phi(s) \Big) ^{\frac{na}{p(n+a)}} \phi(s) ^{1- \frac{na}{(n+a) p}}\\
&\le & \label{eqn:2.16}C_4 \frac{c_\omega^n}{V_\omega} A_{s}^{ \frac{a}{n+a}  } \phi(s) ^{1- \frac{na}{(n+a) p}} + C_5 \frac{c_\omega^n}{V_\omega} A_s \phi(s).
\eea
Here the constant $C_5>0$ depends only on $n, p,\gamma$ and $C_4>0$ depends additionally on $\kappa, K$. Note that the inequality \neweqref{eqn:2.16} holds for {\em any} $s>0$, and the constants $C_4, C_5$ are independent of $s$. We remark that by the choice of $q$, $\frac{na}{(n+a) p } = 1$ when $p\in (0, n)$, and $\frac{na}{(n+a) p } < 1$ when $p=n$, which justifies the H\"older inequality used above.

We now apply Lemma \ref{lemma 2} to conclude that when \begin{equation}\label{eqn:bar s}s\ge \bar s = \max(1, \exp[(2 C_1 C_5 c_\omega^n/V_\omega )^{1/p} ])\end{equation} we have
\begin{equation}\label{eqn:2.17}
\phi(s)\le \frac{C_1}{(\log s)^p} \le \frac 12 \frac{1}{C_5 c_\omega^n/V_\omega}.
\end{equation}
Combining \neweqref{eqn:2.17} and \neweqref{eqn:2.16}, we see that when $s\ge \bar s$
$$
A_s\le 2 C_4 \frac{c_\omega^n}{V_\omega} A_{s}^{\frac{a}{n+a}} \phi(s)^{1-\frac{na}{p(n+a)}}.
$$
Dividing both sides by $A_{s}^{\frac{a}{n+a}}$,  we easily obtain that when $s\ge \bar s$
\begin{equation}\label{eqn:2.18}
A_s\le (2C_4)^{\frac{n+a}{n}} (\frac{c_\omega^n}{V_\omega})^{\frac{n+a}{n}} \phi(s) ^{\frac{n+a}{n} (1- \frac{na}{p(n+a)})}\le C_6(\frac{c_\omega^n}{V_\omega})^{\frac{n+a}{n}}.
\end{equation}
By the definition of $A_s$ in \neweqref{eqn:As}, we see from \neweqref{eqn:2.18} that
\begin{equation}\label{eqn:2.19}
\int_{\Omega_{\bar s}} (-\varphi - \bar s )^{a} e^{nF} \omega_X^n \le C_6 (\frac{c_\omega^n}{V_\omega})^{\frac{a}{n}}.
\end{equation}
Since $a>0$, by the calculus inequality $|x - y|^a\le 2^a(x^a + y^a)$ for any $x,y>0$, we easily obtain from \neweqref{eqn:2.19} that 
\begin{equation}\label{eqn:2.20}
\int_{\Omega_{\bar s}} (-\varphi )^{a} e^{nF} \omega_X^n \le 2^a C_6 (\frac{c_\omega^n}{V_\omega})^{\frac{a}{n}} + C_7 \bar s^a.
\end{equation}
Finally note that on $X\backslash \Omega_{\bar s}$, $0\le -\varphi \le \bar s$, hence  from \neweqref{eqn:2.20}
\bea\nonumber
\int_X (-\varphi)^a e^{nF}\omega_X^n & \le &\int_{\Omega_{\bar s}} (-\varphi)^a e^{nF}\omega_X^n + \int_{X\backslash \Omega_{\bar s}} \bar s^a e^{nF}\omega_X^n\\
&\le \nonumber&2^aC_6 (\frac{c_\omega^n}{V_\omega})^{\frac{a}{n}} + C_8 \bar s^a\\
&=\nonumber &2^aC_6 (\frac{c_\omega^n}{V_\omega})^{\frac{a}{n}} + C_8 \exp[a (2 C_1 C_5 c^n_\omega/V_\omega  )^{1/p}]=:C_e.
\eea
Here $C_e>0$ is the desired constant in \neweqref{eqn:inequality 2} with an explicit dependence on the relative volume $c^n_\omega/V_\omega$.

\medskip

\noindent{\bf Step 4.} We now show the Trudinger-like inequality \neweqref{eqn:inequality}. We take $(\frac{n+a}{n})$-th power on both sides of \neweqref{eqn:2.10} and multiply the resulted inequality by a small constant $\alpha>0$ to be determined. Then it follows that on $\Omega_s$
\begin{equation}\label{eqn:2.21}
\alpha (-\varphi - s)^{\frac{n+a}{n}}\le C_9\alpha A_{s,k}^{1/n} (-\psi_{s,k} + A_{s,k}^{1/a}).
\end{equation}
Taking exponential on both sides of \neweqref{eqn:2.21} and integrating it over $\Omega_s$, we then obtain
\begin{equation}\label{eqn:2.22}
\int_{\Omega_s} e^{\alpha (-\varphi - s)^{\frac{n+a}{n}}}\omega_X^n \le e^{C_9\alpha A_{s,k}^{\frac{n+a}{na}}} \int_{\Omega_s} e^{-C_9 \alpha A_{s,k}^{1/n} \psi_{s,k}}\omega_X^n.
\end{equation} 
Note that by \neweqref{eqn:2.18}, $A_{\bar s} \le C_6 (\frac{c_\omega^n}{V_\omega})^{\frac{n+a}{n}}$ and $A_{\bar s, k}\to A_{\bar s}$, so when $k$ is large enough we have $$A_{\bar s, k} \le 2C_6 (\frac{c^n_\omega}{V_\omega})^{\frac{n+a}{n}}.$$
If we choose $\alpha>0$ small enough such that 
$$C_9 A_{\bar s, k}^{1/n} \alpha \le (2C_6)^{1/n} C_9 (\frac{c_\omega^n}{V_\omega})^{\frac{n+a}{n^2}}\alpha < \alpha(X,\kappa\omega_X)$$
where $\alpha(X,\kappa\omega_X)>0$ is the $\alpha$-invariant of the K\"ahler manifold $(X,\kappa \omega_X)$. Then from \neweqref{eqn:2.22} we get
\begin{equation}\label{eqn:2.24}
\int_{\Omega_{\bar s}} e^{\alpha (-\varphi - \bar s)^{\frac{n+a}{n}}}\omega_X^n \le \exp \Big({C_{10} ( \frac{c_\omega^n}{V_\omega})^{\frac{(n+a)}{n a}} }\Big).
\end{equation}
It is then elementary to see that \neweqref{eqn:2.24} implies 
\begin{equation}\label{eqn:2.25}
\int_{\Omega_{\bar s}} e^{\alpha (-\varphi )^{\frac{n+a}{n}}}\omega_X^n \le \exp \Big({C_{10} ( \frac{c_\omega^n}{V_\omega})^{\frac{(n+a)}{n a}} } + 2 \alpha \bar s^{\frac{n+a}{n}}\Big).
\end{equation}
Again observing that $-\varphi \le \bar s$ on $X\backslash \Omega_{\bar s}$, we conclude from \neweqref{eqn:2.25} that
\begin{equation}\label{eqn:thm 1}
\int_X e^{\alpha (-\varphi )^{\frac{n+a}{n}}}\omega_X^n \le V_{\omega_X} e^{\alpha \bar s^{\frac{n+a}{n}}} + \exp \Big({C_{10} ( \frac{c_\omega^n}{V_\omega})^{\frac{(n+a)}{n a}} } + 2\alpha \bar s^{\frac{n+a}{n}}\Big) =:C_T.
\end{equation}
Since $\bar s$ is explicitly given in \neweqref{eqn:bar s}, the constant $C_T$ has an explicit dependence on $c^n_\omega/V_\omega$. Finally note that $\frac{n+a}{n} =1+ \frac{p}{n-p} = \frac{n}{n-p} = q $. This completes the proof of the inequality \neweqref{eqn:inequality}. Q.E.D.

\section{Proof of Theorem \ref{thm:main3}}
\setcounter{equation}{0}

Again, we drop the subindex $\o$ from $F_\o$ for notational simplicity.
Let $(\varphi, F)$ solve the coupled system \neweqref{eqn:3.1} stated in \S 2. Fix the number $p\in (0,n]$.

\smallskip
\noindent Let $\delta = \frac{1}{10}K_2$. We solve the auxiliary complex Monge-Amp\`ere equation
\begin{equation}\label{eqn:3.3}
(\omega + \ddbar \psi_{k})^n = \frac{\tau_k(-\varphi + \delta F)^p}{A_k} c_\omega^n e^{nF} \omega_X^n,\quad \sup_X\psi_k = 0,
\end{equation}
where 
\begin{equation}\label{eqn:3.4}
A_k = \frac{c_\omega^n}{V_\omega}\int_{X} \tau_k(-\varphi + \delta F)^p e^{nF} \omega_X^n \to \frac{c_\omega^n}{V_\omega}\int_\Omega (-\varphi + \delta F)^p e^{nF} \omega_X^n=:A_\infty,
\end{equation} as $k\to\infty$, 
where $\Omega = \{-\varphi + \delta F >0\}$.
Note that the constant $q>1$ in Theorem \ref{thm:main}, so by H\"older inequality and \neweqref{eqn:inequality 2} in Theorem \ref{thm:main}, we have $\int_X (-\varphi)^p e^{nF}\omega_X^n\le C$ for some constant $C>0$ depending additionally on $K_1$ and $c_\omega^n/V_\omega$. Moreover, by the assumption \neweqref{eqn:assumption},  $\int_\Omega |F |^pe^{nF} \omega_X^n \le K_1$, so we have
$A_\infty \le C$ for some constant depending on $K_1$ and ${c_\omega^n/}{V_\omega}$. Thus $A_k\le C$ for $k$ sufficiently large. 

\medskip

We consider the function $$\Psi = -\varepsilon (-\psi_{k} + \Lambda)^{\frac{n}{n+p}} - \varphi + \delta F.$$
with 
$$\varepsilon = \Big(\frac{(n+p) (1+\delta c_\theta)}{n^2} \Big)^{\frac{n}{n+p}} A_k^{\frac{1}{n+p}} ,$$ and $\Lambda$ is chosen so that $\Lambda^{p/(n+p)} = \frac{2n}{n+p}\varepsilon$, i.e. 
$$\Lambda = (\frac{2n}{n+p})^{\frac{n+p}{p}} \Big(\frac{(n+p) (1+\delta c_\theta)}{n^2} \Big)^{{n/p}} A_k^{1/p}. $$
We claim that $\Psi\le 0$. By the definition of $\Omega$, it suffices to show the case when $\max_X \Psi = \Psi(x_0)$ for some point $x_0\in \Omega$. We calculate at $x_0$ as in the previous section.
\bea\nonumber
0& \ge & G^{i\bar j} \Psi_{i\bar j}\\
&\ge & \frac{\varepsilon n}{n+p} (-\psi_k + \Lambda)^{-\frac{p}{n+p}} G^{i\bar j} (\omega_{\psi_k})_{i\bar j} -  \frac{\varepsilon n}{n+p} (-\psi_k + \Lambda)^{-\frac{p}{n+p}} G^{i\bar j}\omega_{i\bar j}\nonumber\\
&& \nonumber - G^{i\bar j} (\omega_\varphi)_{i\bar j} + G^{i\bar j} \omega_{i\bar j} - \delta c_\theta + \delta G^{i\bar j}\theta_{i\bar j}\nonumber\\
&\ge & \nonumber  \frac{\varepsilon n^2}{n+p} (-\psi_k + \Lambda)^{-\frac{p}{n+p}} \Big( \frac{(-\varphi + \delta F)^p}{A_k}   \Big)^{1/n} - 1 - \delta c_\theta + ( 1- \frac{\varepsilon n \Lambda^{-\frac{p}{n+p}}}{n+p} - \frac{1}{10}  ) G^{i\bar j} \omega_{i\bar j}\\
&\ge & \nonumber  \frac{\varepsilon n^2}{n+p} (-\psi_k + \Lambda)^{-\frac{p}{n+p}} \Big( \frac{(-\varphi + \delta F)^p}{A_k}   \Big)^{1/n} - 1 - \delta c_\theta .
\eea
By the choice of $\varepsilon$, it follows that at $x_0$, 
$$(-\varphi + \delta F)^p \le \Big( \frac{ (n+p) (1+\delta c_\theta)   }{n^2 \varepsilon}   \Big)^n A_k (-\psi_k + \Lambda)^{\frac{np}{n+p}} = \varepsilon^p (-\psi_k + \Lambda)^{\frac{np}{n+p}},$$
that is $\Psi(x_0)\le 0$. This proves the claim that $\Psi\le 0$ on $X$. Since $A_k\le C$ for $k$ large enough, we derive from $\Psi\le 0$ that in $X$
$$\delta F \le -\varphi + \delta F\le C(-\psi_k + 1)^{\frac{n}{n+p}}\le - \epsilon \psi_{k} + C_\epsilon, $$
where the last inequality follows from the elementary inequality $(x+1)^{n/(n+p)}\le \epsilon x + C_\epsilon$ for any $\epsilon>0$. In particular, this shows that for {\em any} $r\ge 1$
\begin{equation}\label{eqn:iteration}\int_X e^{r \delta F}\omega_X^n \le C \int_X e^{-r \epsilon \psi_k} \omega_X^n\le C_r,\end{equation}
where we choose $\epsilon>0$ small so that $\epsilon r<$ the alpha-invariant of $(X,\kappa \omega_X)$. In particular, this implies that $e^{nF}$ is bounded in  $L^{p'}(X,\omega_X^n)$ for any $p'> 1$. An application of \cite{GPT} then implies that the $L^\infty$ norm of $\varphi$ is bounded by a constant depending only $\| e^{nF}\|_{L^{^{p'}}(X,\omega_X^n)}$ for some $p'>1$ and $c_\omega^n/V_\omega$ (c.f. Theorem \ref{thm:main2} in \S 2), and subsequently by \neweqref{eqn:iteration}, $\| \varphi\|_{L^\infty}$ depends only on $\kappa, K_1, K_2,c_\theta,\gamma$ and $c_\omega^n/V_\omega$.

\medskip

Moreover, under the assumption \neweqref{eqn:assumption}, we will show $\sup_X F$ is bounded above by a uniform constant. To see this, we begin with a mean-value type inequality which was proved in \cite{GPS} for complex Monge-Amp\`ere equations and the arguments there can be easily adapted to the current situation. But for convenience of the readers, we include a proof of this inequality in Lemma \ref{lemma:GPS} below. 

\begin{lemma}\label{lemma:GPS}
Let $u\in C^2(X)$ be a $C^2$ function on $X$ that satisfies the differential inequality 
\begin{equation}\label{eqn:diff}\Box_{\omega_\varphi} u \ge - a ,\end{equation}
then there is a constant $C>0$ depending on $n,\kappa, \gamma, \kappa, p, K_1, c_\omega^n/V_\omega$ and $a\ge 0$ such that 
$$\sup\nolimits_X u \le C( 1 + \frac{c_\omega^n}{V_\omega} \int_X |u| e^{nF} \omega_X^n ).$$
\end{lemma}

\noindent {\em Proof of Lemma \ref{lemma:GPS}}. 
As in \cite{GPS}, we may assume that $N: = \frac{c_\omega^n}{V_\omega} \int_X |u| e^{nF} \omega_X^n\le 1$, otherwise, we can consider the rescaled function $\tilde u = u/N$, which still satisfies the differential {\em inequality} \neweqref{eqn:diff} with the {\em same} $a$. We also assume $\{u>0\}\neq \emptyset$, otherwise this lemma is trivial.

Let $s>0$ be a positive number such that the super-level set $U_s: = \{u>s\}$ is non-empty. 
We consider the auxiliary Monge-Amp\`ere equation
\begin{equation}\label{eqn:3.6}
(\omega + \ddbar \psi_{s,k})^n = \frac{\tau_k(u-s)}{A_{s,k}} c_\omega^n e^{nF} \omega_X^n,\quad \sup_X\psi_{s,k} = 0,
\end{equation}
where as $k\to\infty$
$$A_{s,k} = \frac{c_\omega^n}{V_\omega} \int_X \tau_k(u-s) e^{nF} \omega_X^n\to A_{s} = \frac{c_\omega^n}{V_\omega} \int_{U_s} (u-s) e^{nF} \omega_X^n.$$ The condition that $N\le 1$ implies that $A_{s}\le 1$, so if $k>1$ is sufficiently large, we have $A_{s,k}\le 2$. 

The following argument is  similar to that in {\bf Step 2} in the proof of Theorem \ref{thm:main}. We have known that $\sup_X |\varphi|$ is bounded. So take $\Lambda_0 = \sup_X|\varphi| + 1$. We consider the  function $$\Phi_u: = - \varepsilon (-\psi_{s,k} + \varphi + \Lambda_0)^{\frac{n}{n+1}} + u - s,$$
and we claim that for $\varepsilon =\varepsilon(s,k,a)>0 $ satisfying the equation
\begin{equation}\label{eqn:3.7}\varepsilon^{n+1} = A_{s,k} ( a + \frac{n\varepsilon}{n+1}  )^n,\end{equation} we have $\sup_X \Phi_u\le 0$. First we observe that from equation \neweqref{eqn:3.7} and $A_{s,k}\le 2$ it holds that $\varepsilon \le C A_{s,k}^{1/(n+1)}$ for some $C = C(n, a)>0$. 
If the maximum of $\Phi_u$ is achieved at some point $x_0\in U_s$, then at this point by maximum principle
\bea\nonumber
0& \ge &\Box_{\omega_\varphi} \Phi_u(x_0)\\
&\ge \nonumber& \frac{n\varepsilon}{n+1} (-\psi_{s,k} + \varphi + \Lambda_0) ^{-\frac{1}{n+1}} \Big( G^{i\bar j} (\omega_{\psi_{s,k}})_{i\bar j} - G^{i\bar j} (\omega_\varphi)_{i\bar j}   \Big) - a\\
&\ge \nonumber& \frac{n\varepsilon}{n+1} (-\psi_{s,k} + \varphi + \Lambda_0) ^{-\frac{1}{n+1}} \Big( n   ( \frac{u - s }{A_{s,k}}  )^{1/n} - 1         \Big) - a  \\
&\ge \nonumber& \frac{n^2\varepsilon}{n+1} (-\psi_{s,k} + \varphi + \Lambda_0) ^{-\frac{1}{n+1}}  \Big( \frac{u - s }{A_{s,k}} \Big )^{1/n}  - \frac{n\varepsilon}{n+1}      - a  .
\eea
It is then elementary to see that $\Phi_u\le 0$ by the choice of $\varepsilon$ in \neweqref{eqn:3.7}. This immediately implies
$$\int_{U_s} \exp\Big( \alpha \frac{(u-s)^{\frac{n+1}{n}}}{A_{s,k}^{1/n}} \Big) \omega_X^n \le C(n,a) \int_{U_s} e^{ -\alpha C(n,a) \psi_{s,k} } \omega_X^n \le C, $$
if we choose $\alpha>0$ small enough so that $\alpha C(n,a)$ is less  than the alpha invariant of $(X,\kappa \omega_X)$. Letting $k\to \infty$ we obtain 
\begin{equation}\label{eqn:3.8}
\int_{U_s} \exp\Big( \alpha \frac{(u-s)^{\frac{n+1}{n}}}{A_{s}^{1/n}} \Big) \omega_X^n \le C.
\end{equation}
This equation together with H\"older-Young's inequality yield that for any $r>n$
\begin{equation}\label{eqn:3.9}
\int_{U_s} (u-s)^{\frac{(n+1)r}{n}} e^{nF} \omega_X^n \le C A_s^{r/n}.
\end{equation}
From now on we fix an $r >n$.
Then we can apply H\"older inequality to obtain 
\bea
A_s & \nonumber= & \frac{c^n_\omega}{V_\omega}\int_{U_s} (u-s) e^{nF} \omega_X^n \\
& \le  \nonumber& \frac{c_\omega^n}{V_\omega}\Big ( \int_{U_s} (u-s)^{\frac{(n+1)r}{n}} e^{nF} \omega_X^n\Big)^{\frac{n}{r (n+1)}} \Big(\int_{U_s} e^{nF} \omega_X^n \Big)^{1- \frac{n}{r (n+1)}}\\
&\le & \nonumber\frac{c_\omega^n}{V_\omega} C A_s^{1/(n+1)} \phi_u(s) ^{1- \frac{n}{r (n+1)}},
\eea
where $\phi_u(s) = \int_{U_s} e^{nF} \omega_X^n$. Then we have
$$A_s\le C (\frac{c_\omega^n}{V_\omega})^{1 +\frac 1n} \phi_u(s)^{1 + \frac{r-n}{rn}}.  $$
This combined with the definition of $A_s$ implies that for any $t>0$, there exists a constant $\bar C>0$ depending on $n, \gamma, \kappa, a, K_1, K_2,$ such that 
$$ t \phi_u(s+t) \le \bar C (\frac{c_\omega^n}{V_\omega})^{1/n} \phi_u(s)^{1+ \frac{r-n}{rn}}.$$
Let $s_0>0$ be a number such that $\bar C (\frac{c_\omega^n}{V_\omega})^{1/n} \phi_u(s_0)^{\frac{r-n}{rn}}\le 1/2$. This $s_0$ can be chosen as $(2\bar C)^{\frac{rn}{r-n}} (\frac{c_\omega^n}{V_\omega})^{n/(r-n)}$, since by the assumption $N\le 1$
$$\phi_u(s_0) \le \frac{N}{s_0} \frac{1}{c_\omega^n/V_\omega}\le  \frac{1}{s_0} \frac{1}{c_\omega^n/V_\omega}.$$ Then a De Giorgi type iteration argument of Kolodziej \cite{K} (see also \cite{GPT} and \cite{Di}) implies that $\phi_u(s) = 0$ for $s\ge S_\infty$ for some uniform constant $$S_\infty = s_0 + \frac{1}{1- 2^{-{(r-n)/}{rn }}},$$ and this gives $u\le S_\infty$ as desired. Q.E.D. 

\medskip

We now apply Lemma \ref{lemma:GPS} to $u: = F - K_2 \varphi$, which satisfies
$$\Box_{\omega_\varphi} u = -c_\theta + G^{i\bar j}\theta_{i\bar j} - K_2 + K_2 G^{i\bar j} \omega_{i\bar j}\ge -c_\theta - K_2 = : -  a. $$
Lemma \ref{lemma:GPS} yield the upper bound of $u$ (hence that of $F$) in terms of its $L^1$ integral, while the latter is bounded since $e^{nF}$ is bounded in $L^r(\omega_X^n)$ for any $r>1$ and $\varphi$ is bounded in $L^\infty$ by \cite{GPT}.

\medskip

If in addition we assume \neweqref{eqn:K3}, that is
$$\theta \le K_3 \omega,$$
then apply Lemma \ref{lemma:GPS} to $u: = - F -  K_3 \varphi$, we can also get a uniform {\em lower} bound of $F$ depending on $K_3$. Q.E.D.

\bigskip

\noindent Department of Mathematics \& Computer Science, Rutgers University, Newark, NJ 07102 USA

\noindent bguo@rutgers.edu,

\medskip

\noindent Department of Mathematics, Columbia University, New York, NY 10027 USA

\noindent phong@math.columbia.edu

\end{document}